\documentclass[10.5pt]{amsart}
\usepackage[english]{babel}
\usepackage{amsmath}
\usepackage{amssymb}
\usepackage{amsthm}
 \newtheorem{thm}{Theorem}[section]
 \newtheorem{cor}[thm]{Corollary}
 \newtheorem{lem}[thm]{Lemma}
 \newtheorem{prop}[thm]{Proposition}
 \theoremstyle{definition}
 \newtheorem{defn}[thm]{Definition}
 \theoremstyle{remark}
 \newtheorem{rem}[thm]{Remark}
 \newtheorem*{ex}{Example}
 \numberwithin{equation}{section}

\newcommand{\scal}[1]{\left<#1\right>}
\newcommand{\Hq}{\mathbb H}
\newcommand{\Sq}{\mathbb S}

\newcommand{\N}{\mathbb{N}}

\newcommand{\R}{\mathbb{R}}      

\newcommand{\C}{\mathbb{C}}

\newcommand{\FN}{\mathcal{F}^N_{Slice}(\Hq)}

\newcommand{\FNi}{\mathcal{F}^{N}_{I}(\Hq)}
\newcommand{\FNj}{\mathcal{F}^{N}_{J}(\Hq)}
\newcommand{\FNjj}{\mathcal{F}^{N}_{-J}(\Hq)}
\newcommand{\AN}{\mathcal{A}^N_{Slice}(\mathbb{B})}

\title[]{On slice polyanalytic functions of a quaternionic variable}

\author[D. Alpay]{Daniel Alpay}
\address{(DA) Schmid College of Science and Technology, Chapman University, Orange 92866, CA, US}
\email{alpay@chapman.edu}

\author[K. Diki]{Kamal Diki}
\address{(KD) Marie Sklodowska-Curie fellow of the Istituto Nazionale di Alta Matematica \\ Politecnico di
Milano\\Dipartimento di Matematica\\Via E. Bonardi, 9\\20133 Milano,
Italy}
\email{kamal.diki@polimi.it}

\author[I. Sabadini]{Irene Sabadini}
\address{(IS) Politecnico di
Milano\\Dipartimento di Matematica\\Via E. Bonardi, 9\\20133 Milano\\Italy}
\email{irene.sabadini@polimi.it}

\begin{document}
\maketitle
\begin{abstract}

In this paper, we introduce the quaternionic slice polyanalytic functions and we prove some of their properties. Then, we apply the obtained results to begin the study of the quaternionic Fock and Bergman spaces in this new setting. In particular, we give explicit expressions of their reproducing kernels.

 \end{abstract}



\section{Introduction}
The theory of polyanalytic functions is an interesting topic in complex analysis. It extends the concept of holomorphic functions to nullsolutions of higher order powers of the Cauchy-Riemann operator. An excellent reference on this subject is the book of Balk \cite{Balk1991}. Some famous Hilbert spaces of holomorphic functions that were extended to the setting of polyanalytic functions are the Bergman and Fock spaces, see for example \cite{AF2014,Alpay2015,Balk1991,Koselev1977} and the references therein. In the last years, the classical theory of holomorphic functions in complex analysis was extended to obtain the theory of slice hyperholomorphic functions of a quaternionic variable, see \cite{ColomboSabadiniStruppa2011,GentiliStoppatoStruppa2013}  and the references therein. This function theory has several applications, in particular in operator theory and Schur analysis, see \cite{ACS2016,ColomboSabadiniStruppa2011}. In this paper, we extend the definition of slice hyperholomorphic functions to higher order and define the slice polyanalytic functions of a quaternionic variable. Then, we shall use the obtained results to introduce and study the Fock and Bergman spaces of quaternionic slice polyanalytic functions and give explicit formulas for their reproducing kernels. Note that by considering polyanalytic functions with respect to the classical Cauchy-Fueter regularity on quaternions, it turns out that even the simple example given by $F(q,\overline{q})=|q|^2$ is not polyanalytic of order 2. However, a natural question that may arise in this direction is about a possible extension of the well-known Fueter mapping theorem on quaternions allowing to construct Cauchy-Fueter polyanalytic functions starting from slice polyanalytic functions of the same order. The paper has the following structure: in  Section 2, we review some useful preliminaries on classical complex polyanalytic functions and on quaternionic slice hyperholomorphic functions. Then, in Section 3 we introduce the quaternionic slice polyanalytic functions. In particular, on slice domains we show that a slice polyanalytic function is the sum of the quaternionic conjugate powers multiplied by slice regular functions, thus extending the analogous result for complex functions. We prove also the counterparts of the Splitting Lemma, Identity Principle, Representation Formula, Extension Lemma and the Refined Splitting Lemma in this framework. We also discuss  slice polyanalytic functions as a subclass of slice functions on axially symmetric domains. In particular, we prove a version of the identity principle in this situation. In Section 4, we introduce and study the Fock space of slice polyanalytic functions on quaternions and we give the formula of its reproducing kernel. Section 5 treats the case of the Bergman theory of the second kind in the quaternionic slice polyanalytic setting  in the case of the unit ball. We conclude the paper by a brief discussion related to avenues for further research.

\section{Preliminaries}
 In this section, we revise the needed material concerning complex polyanalytic functions and the theory of slice regular functions on quaternions.
 \subsection{Polyanalytic functions in classical complex analysis}
 We begin by providing the background on complex polyanalytic functions.
The reader interested in more details, may consult the book \cite{Balk1991}.
\begin{defn}
Let $\Omega$ be a domain of $\C$. A function $f:\Omega\longrightarrow \C$ is said to be polyanalytic of order $n$ or $n-$analytic if $$\displaystyle\left(\frac{\partial}{\partial \overline{z}}\right)^nf(z)=0, \text{ } \forall z\in\Omega.$$
The space of all polyanalytic functions of order $n$ is denoted by $H_n(\Omega).$
\end{defn}
\begin{ex}
The function $F(z)=1-z\overline{z}$ is polyanalytic of order $2$ on $\C$.
\end{ex}
\begin{prop} \label{poly}
Let $\Omega$ be a domain of $\C$ and $f:\Omega\longrightarrow \C$. Then, the two following conditions are equivalent
\begin{enumerate}
\item $f$ is polyanalytic of order $n$.
\item $\displaystyle f(z)=\sum_{k=0}^{n-1}\overline{z}^ka_k(z), \forall z\in\Omega$ where $a_0,...,a_{n-1}$ are analytic on $\Omega$.
\end{enumerate}
\end{prop}
\begin{prop}
Let $f$ and $g$ be two polyanalytic functions of order $n$ on a domain $\Omega$. If $\Omega_0$ is a subdomain of $\Omega$ such that $f$ and $g$ coincide on $\Omega_0$, then $f$ and $g$ coincide everywhere in $\Omega$.
\end{prop}
In the book of Balk \cite{Balk1991}, the Fock space $\mathcal{F}_n(\C)$ of polyanalytic functions of order $n$ is defined by $$\mathcal{F}_n(\C)=\lbrace{f\in H_n(\C);\int_\C\vert{f(z)}\vert^2e^{-\vert{z}\vert^2}d\lambda(z)<\infty}\rbrace, $$
where $d\lambda(z)$ denotes the usual Lebesgue measure on the complex plane. Note that, $\mathcal{F}_n(\C)$ is a reproducing kernel Hilbert space whose reproducing kernel is
\begin{equation}
\displaystyle F_n(z,w)=e^{\overline{w}z}\sum_{k=0}^{n-1}(-1)^k{n \choose k+1}\frac{1}{k!}\vert{z-w}\vert^{2k}
\end{equation}
Moreover, for all $f\in\mathcal{F}_n(\C)$ and $z\in\C$ we have $$\vert{f(z)}\vert\leq \sqrt{n}e^{\frac{\vert{z}\vert^2}{2}}\Vert{f}\Vert_{\mathcal{F}_n(\C)}.$$
On the other hand, the Bergman space $A^2_n(\mathbb{D})$ of polyanalytic functions of order $n$ in the unit disc is given by $$\displaystyle A^2_n(\mathbb{D})=\lbrace{f\in H_n(\mathbb{D});\int_{\mathbb{D}}\vert{f(z)}\vert^2d\lambda(z)<\infty}\rbrace.$$
Also
$A^2_n(\mathbb{D})$ is a reproducing kernel Hilbert space whose reproducing kernel is given by
\begin{equation}
\displaystyle B_n(z,w)=\frac{n}{\pi(1-\overline{w}z)^{2n}}\sum_{k=0}^{n-1}(-1)^k{n \choose k+1}{n+k \choose n}\vert{1-\overline{w}z}\vert^{2(n-1-k)}\vert{z-w}\vert^{2k}
\end{equation}
for any $z,w\in\mathbb{D}$.
Moreover, for all $f\in A^2_n(\mathbb{D})$ and $z\in\mathbb{D}$, we have the following
$$\vert{f(z)}\vert\leq \frac{n}{\sqrt{\pi}}\frac{\Vert{f}\Vert_{A^2_n(\mathbb{D})}}{(1-\vert{z}\vert^2)}.$$
 \subsection{Quaternions and slice regular functions}
For more details about the theory of slice regular functions, originally introduced in \cite{GentiliStruppa07} and its different generalizations and applications one can see for example \cite{ColomboSabadiniStruppa2016,ColomboSabadiniStruppa2011,
GentiliStoppatoStruppa2013} and the references therein. In this section we summarize only the results needed in the sequel.\\ \\
The non-commutative field of quaternions is defined to be
$$\Hq=\lbrace{q=x_0+x_1i+x_2j+x_3k\quad ; \ x_0,x_1,x_2,x_3\in\R}\rbrace$$ where the imaginary units satisfy the multiplication rules $$i^2=j^2=k^2=-1\quad \text{and}\quad ij=-ji=k, jk=-kj=i, ki=-ik=j.$$
On $\Hq$ the conjugate and the modulus of $q$ are defined respectively by
$$\overline{q}=Re(q)-Im(q) \quad \text{where} \quad Re(q)=x_0, \quad Im(q)=x_1i+x_2j+x_3k$$
and $$\vert{q}\vert=\sqrt{q\overline{q}}=\sqrt{x_0^2+x_1^2+x_2^2+x_3^2}.$$
Note that the quaternionic conjugation satisfy the property $\overline{ pq }= \overline{q}\, \overline{p}$ for any $p,q\in \Hq$.
Moreover, the unit sphere $$\lbrace{q=x_1i+x_2j+x_3k;\text{ } x_1^2+x_2^2+x_3^2=1}\rbrace$$ coincides with the set of all  imaginary units given by $$\mathbb{S}=\lbrace{q\in{\Hq};q^2=-1}\rbrace.$$
  Any quaternion $q\in \Hq\setminus \R$ can be written in a unique way as $q=x+I y$ for some real numbers $x$ and $y>0$, and imaginary unit $I\in \mathbb{S}$, in fact $q=x_0+\dfrac{x_1i+x_2j+x_3k}{|x_1i+x_2j+x_3k|}|x_1i+x_2j+x_3k|$.
Then, for every given $I\in{\mathbb{S}}$, the slice $\C_I$ is defined to be $\mathbb{R}+\mathbb{R}I$ and it is isomorphic to the complex plane $\C$ so that it can be considered as a complex plane in $\Hq$ passing through $0$, $1$ and $I$. The semi-slice $\C_I^+$ is given by the set $\lbrace{x+yI;x,y\in\R, y\geq 0}\rbrace$. If $q=x_0\in\mathbb R$ then $q\in\mathbb C_I$ for all $I\in\mathbb S$. It is immediate that  $\Hq=\underset{I\in{\mathbb{S}}}{\cup}\C_I.$
\\ In \cite{GentiliStruppa07}, the authors proposed a new definition to extend the classical theory of holomorphic functions in complex analysis to the quaternionic setting. This leads to the new theory of slice hyperholomorphic or slice regular functions on quaternions:
\begin{defn}
A real differentiable function $f: \Omega \longrightarrow \Hq$, on a given domain $\Omega\subset \Hq$, is said to be a (left) slice  regular function if, for every $I\in \Sq$, the restriction $f_I$ to the slice $\C_{I}$ satisfies
$$
\overline{\partial_I} f(x+Iy):=
\dfrac{1}{2}\left(\frac{\partial }{\partial x}+I\frac{\partial }{\partial y}\right)f_I(x+Iy)=0,
$$
on $\Omega_I$. The slice derivative $\partial_S f$ of $f$ is defined by :
\begin{equation*}
\partial_S(f)(q):=
\left\{
\begin{array}{rl}
\partial_I(f)(q)& \text{if } q=x+Iy, y\neq 0\\
\displaystyle\frac{\partial}{\partial{x}}(f)(x) & \text{if } q=x \text{ is real}.
\end{array}
\right.
\end{equation*}
\end{defn}
In addition we introduce the following terminology
\begin{defn}
\begin{enumerate}
\item A quaternionic valued function on a domain $\Omega$ is said to be (quaternionic) intrinsic if $f(\Omega_I)\subset\C_I$ for any $I\in\mathbb{S}$. We note that if a function is expressed by a power series or, in particular, by a polynomial then it is intrinsic if and only if it has real coefficients.
\item  A function which is slice regular on the whole space of quaternions $\Hq$ is said to be entire.
\end{enumerate}
\end{defn}
We will refer to left slice regular functions as slice regular functions, for short. The set of these functions is denoted by $\mathcal{SR}(\Omega)$. It turns out that $\mathcal{S}\mathcal{R}(\Omega)$ is a right vector space over the noncommutative field $\Hq$.
\begin{rem}
The multiplication and composition of slice regular functions are not slice regular, in general. Moreover, the slice derivative does not satisfy the Leibniz rule with respect to the pointwise multiplication. However, the composition $f\cdot g$ of two slice regular functions is slice regular if $g$ is intrinsic and the pointwise product $fg$ is slice regular if $f$ is intrinsic, see \cite{ColomboSabadiniStruppa2011}.
\end{rem}

 According to the definition, the basic polynomials in $q$ with quaternionic coefficients on the right are slice regular. Moreover, for any power series $\displaystyle\sum_n q^na_n$, there exists $0\leq R \leq \infty$, called the radius of convergence such that the power series is a slice regular function on $B(0,R):= \{q\in \Hq; \, |q|<R\}$. The space of slice regular functions is endowed with the natural topology of uniform convergence on compact sets. The characterization of slice regular functions on a ball $B = B(0,R)$ centered at the origin is given by

\begin{thm}[Series expansion]
An $\Hq$-valued function $f$ is slice regular on $B(0,R)\subset \Hq$ if and only if it has a series expansion of the form:
$$f(q)=\sum_{n=0}^{+\infty} q^n\frac{1}{n!}\partial^{(n)}_S(f)(0)$$
converging on $B(0,R)=\{q\in\Hq;\mid q\mid<R\}$.
\end{thm}

\begin{defn}
A domain $\Omega\subset \Hq$ is said to be a slice domain (or just $s$-domain) if  $\Omega\cap{\mathbb{R}}$ is nonempty and for all $I\in{\mathbb{S}}$, the set $\Omega_I:=\Omega\cap{\C_I}$ is a domain of the complex plane $\C_I$.
If moreover, for every $q=x+Iy\in{\Omega}$, the whole sphere $x+y\mathbb{S}:=\lbrace{x+Jy; \, J\in{\mathbb{S}}}\rbrace$
is contained in $\Omega$, we say that  $\Omega$ is an axially symmetric slice domain.
\end{defn}

\begin{ex}
The whole space $\Hq$  and the Euclidean ball $B=B(0,R)$ of radius $R$ centered at the origin are axially symmetric slice domains.
\end{ex}
The following properties of slice regular functions are of particular interest and will be very useful for the next sections of this paper, see \cite{ColomboSabadiniStruppa2011,GentiliStoppatoStruppa2013}.

\begin{lem}[Splitting Lemma]\label{split} Let $f$ be a slice regular function on a domain $\Omega$. Then, for every $I,J\in\mathbb S$ with $I\perp J$, there exist two holomorphic functions $F,G:\Omega_{I}\longrightarrow{\C_I}$ such that for all $z=x+Iy\in{\Omega_I}$, we have
$$f_I(z)=F(z)+G(z)J,$$
 where $\Omega_I=\Omega\cap{\C_I}$ and $\C_I=\mathbb{R}+\mathbb{R}I.$
\end{lem}

\begin{thm}[Representation Formula]\label{repform}
Let $\Omega$ be an axially symmetric slice domain and $f\in{\mathcal{SR}(\Omega)}$. Then, for any $I,J\in{\mathbb{S}}$, we have the formula
$$
f(x+Jy)= \frac{1}{2}(1-JI)f_I(x+Iy)+ \frac{1}{2}(1+JI)f_I(x-Iy)
$$
for all $q=x+Jy\in{\Omega}$.
\end{thm}

\begin{thm}[Identity Principle] \label{IdentityPrin}
Let $f$ and $g$ be two slice regular functions on a slice
domain $\Omega$. If, for some $I\in\mathbb{S}$, $f$ and $g$ coincide on a subset of $\Omega_I$ having an accumulation point in $\Omega_I$, then $f=g$ on the whole domain $\Omega$.
\end{thm}
\begin{lem}[Extension Lemma]\label{extensionLem}
Let $\Omega_I$ be a domain in $\C_I$ symmetric with respect to the real axis and such that $\Omega_I\cap\R\neq \emptyset$. Let $h:\Omega_I\longrightarrow \Hq$ be a holomorphic function. Then, the function $ext(h)$ defined by
$$ext(h)(x+Jy):= \dfrac{1}{2}[h(x+Iy)+h(x-Iy)]+\frac{JI}{2}[h(x-Iy)-h(x+Iy)];  \quad J\in \mathbb{S},$$
extends $h$ to a slice regular function $ext(h)$ on $\overset{\sim}\Omega=\underset{x+Iy\ ; \ x+Jy\in{\Omega}}\cup x+Iy$, the symmetric completion of $\Omega_I$.
Moreover, $ext(h)$ is the unique slice regular extension of $h$.
\end{lem}
In \cite{AlpayColomboSabadini2014}, the authors defined the slice hyperholomorphic Fock space on quaternions as
$$\mathcal{F}^{2}_{Slice}(\Hq):=\lbrace{f\in{\mathcal{SR}(\Hq); \,  \frac{1}{\pi}\int_{\C_I}\vert{f_I(q)}\vert^2 e^{-\vert{q}\vert^2}d\lambda_I(q) <\infty}}\rbrace,$$
 where $f_I = f|_{\C_I}$ and $d\lambda_I(q)=dxdy$ for $q=x+Iy$.
It was proved that $\mathcal{F}^{2}_{Slice}(\Hq)$ is a right quaternionic reproducing kernel Hilbert space whose reproducing kernel is given by \begin{equation}\label{expRK}
e_*(p\overline{q}) = \sum_{n=0}^\infty\frac{p^n\overline{q}^n}{n!}.
\end{equation}
 Equivalently, the reproducing kernel of the slice hyperholomorphic Fock space could be obtained also by taking the slice regular extension of the complex function $e^{z\overline{q}}$ where $z$ and $q$ belong to the same slice. This means that
 \begin{equation}
 e_*(p\overline{q})=ext(e^{z\overline{q}})(p).
 \end{equation}
 The slice Bergman space of the first and second kind were introduced in  \cite{ColomboCervantesSabadini2015,ColomboCervantesLunaSabadiniShapiro2013}. In this paper, we focus on the case
 of the Bergman space of the second kind of the quaternionic unit ball $\mathbb{B}$. For $I\in\mathbb{\Sq}$, the slice hyperholomorphic Bergman space of the second kind is defined to be $$\displaystyle\mathcal{A}_{Slice}(\mathbb{B}):=\lbrace{f\in\mathcal{SR}(\mathbb{B});\int_{\mathbb{B}_I}\vert{f_I(p)}\vert^2d\lambda_I(p)<\infty}\rbrace.$$
 Note that, $\mathcal{A}_{Slice}(\mathbb{B})$ is a right quaternionic Hilbert space which does not depend on the choice of the imaginary unit $I$. Its reproducing kernel is obtained by extending the classical kernel in complex analysis; in closed form it can be written as follows, see \cite{ColomboCervantesSabadini2015}:
 \begin{equation}
 \displaystyle B_S(q,r)=\frac{1}{\pi}(1-2\bar{q}\bar{r}+\bar{q}^2\bar{r}^2)(1-2Re(q)\bar{r}+\vert{q}\vert^2\bar{r}^2)^{-2}.
 \end{equation}
 We note that this kernel can be written also in the following form
\begin{equation}
 \displaystyle B_S(q,r)=\frac{1}{\pi}(1-2Re(r)q+\vert{r}\vert^2q^2)^{-2}(1-2qr+q^2r^2).
 \end{equation}
 \section{Quaternionic slice polyanalytic functions}
In this section, we extend to higher order the theory of slice regular functions on quaternions. First, we start by considering the following simple example
\begin{ex}
 For any $q\in\Hq$, let $F(q)=1-\overline{q}qj$. Then, we have  $$\overline{\partial}_IF_I(x+Iy)=-(x+Iy)j\text{ and }\overline{\partial}^2_IF_I(x+Iy)=0;\text{ }\forall I\in\Sq.$$ We say that $F$ is slice polyanalytic of order $2$ on $\Hq$.
\end{ex}
 The last example, suggests to consider this more general definition :
\begin{defn}
Let $\Omega$ be an open set in $\Hq$ and let $f:\Omega\longrightarrow \Hq$ be a real function $N$ times differentiable. For each $I\in\Sq$, let $\Omega_I=\Omega\cap\C_I$ and let $f_I=f_{\vert_{\Omega_I}}$ be the restriction of $f$ to $\Omega_I$. The restriction $f_I$ is called (left) polyanalytic of order $N$ if it satisfies on $\Omega_I$ the equation $$
\overline{\partial_I}^N f(x+Iy):=
\frac{1}{2^N}\left(\frac{\partial }{\partial x}+I\frac{\partial }{\partial y}\right)^Nf_I(x+Iy)=0.
$$
The function $f$ is called left slice polyanalytic of order $N$, if for all $I\in\Sq$, $f_I$ is left polyanalytic of order $N$ on $\Omega_I$.
\end{defn}
\begin{rem} We note that when dealing with left slice polyanalytic functions we will refer to them simply as slice polyanalytic.
 Due to the lack of commutativity on $\Hq$, we can define in an analogous way the right slice polyanalytic functions on quaternions.
\end{rem}
The set of all slice polyanalytic functions of order $N$ on a domain $\Omega$ is a right vector space over the noncommutative field of quaternions. It will be denoted $\mathcal{SP}_N(\Omega)$ or simply $\mathcal{SP}(\Omega)$ if no confusion can arise with respect to the order. A simple observation that will be needed in the sequel is the following
\begin{prop} \label{Mult}
If $f$ is an intrinsic, slice polyanalytic function of order $N$ and $g$ is a slice regular function on a domain $\Omega$ then the pointwise multiplication $h(q)=f(q)g(q)$ defines also a slice polyanalytic function of order $N$ on $\Omega$.
\end{prop}
\begin{proof}
This holds because $f$ is intrinsic, thus we can use the Leibniz rule.
\end{proof}
\begin{prop}[Splitting Lemma]\label{split} Let $f$ be a slice polyanalytic function of order $N$ on a domain $\Omega\subseteq\Hq$. Then, for any imaginary units $I$ and $J$ with $I\perp J$ there exist $F,G:\Omega_{I}\longrightarrow{\C_I}$ polyanalytic functions of order $N$ such that for all $z=x+Iy\in\Omega_I$, we have
$$f_I(z)=F(z)+G(z)J.$$
 \end{prop}
\begin{proof}
Let $I,J\in\Sq$ be such that $I\perp J$, then $\lbrace{1,I,J,IJ}\rbrace$ forms an orthogonal basis of $\Hq$. Hence, for any $z=x+Iy$ we can write $$f_I(z)=f_0(z)+f_1(z)I+f_2(z)J+f_3(z)IJ$$ where $f_0,..,f_3$ are real valued. This leads to $$f_I(z)=F(z)+G(z)J$$ with $F(z)=f_0(z)+f_1(z)I$ and $G(z)=f_2(z)+f_3(z)I$. However, $f$ is slice polyanalytic of order $N$ which means that $\overline{\partial}_I^Nf_I(x+Iy)=0$ on $\Omega_I$. Thus, by linearity of the operator $\displaystyle\overline{\partial}_I^N$ and linear independence of the basis elements we have $\overline{\partial}_I^NF(x+Iy)=0$ and $\overline{\partial}_I^NG(x+Iy)=0$ on $\Omega_I$. This ends the proof.
\end{proof}
\begin{prop}\label{carac1}
Let $f_0,...,f_{N-1}$ be slice regular functions on a domain $\Omega\subseteq\mathbb H$. Then, the function defined by \begin{equation}\label{sum}
f(q):=\displaystyle\sum_{k=0}^{N-1}\overline{q}^kf_k(q)
\end{equation}
is slice polyanalytic of order $N$ on $\Omega$.
\end{prop}
\begin{proof}
Let $I\in\Sq$ and choose $J\in\Sq$ with $I\perp J$. The Splitting Lemma for slice regular functions yields $${f_{k}}_{\vert_{\C_I}}(x+Iy)=F_k(x+Iy)+G_k(x+Iy)J;\text{ } \forall k=0,...,N$$
where $F_k$ and $G_k$ are $\C_I$ valued holomorphic functions on $\Omega_I$. Hence, we have \begin{align*}
f_I(x+Iy)
 &=\displaystyle \sum_{k=0}^{N-1}(x-Iy)^k{f_{k}}_{\vert_{\C_I}}(x+Iy)
 \\& = \sum_{k=0}^{N-1}(x-Iy)^kF_k(x+Iy)+\sum_{k=0}^{N-1}(x-Iy)^kG_k(x+Iy) J
 \\& = F(x+Iy)+G(x+Iy)J.
 \end{align*}
  It is immediate that $F$ and $G$ are polyanalytic of order $N$ on $\Omega_I$. Thus $\overline{\partial}_I^Nf_I(x+Iy)=0$ on $\Omega_I$ for any $I\in\Sq$.
\end{proof}
Conversely, we have the following
\begin{prop}\label{carac2}
If $f$ is a slice polyanalytic function of order $N$ defined on a slice domain $\Omega\subset\Hq$. Then,
\begin{equation} \label{sum'}
f(q)=\displaystyle\sum_{k=0}^{N-1}\overline{q}^kf_k(q)
\end{equation}
where $f_0,...,f_{N-1}$ are slice regular functions on $\Omega$.
\end{prop}
\begin{proof}
Let $f:\Omega\longrightarrow \Hq$ be a slice polyanalytic function of order $N$ on a slice domain $\Omega$. Then, we can use the Splitting Lemma (Proposition \ref{split}) to write $f_I(z)=F(z)+G(z)J$, where $z=x+yI$ and $F,G$ are $\C_I-$valued polyanalytic functions of order $N$ on $\Omega_I$. Then, by classical complex analysis we have that $$\displaystyle F(z)=\sum_{k=0}^{N-1}\overline{z}^k\varphi_k(z)\text{ and } G(z)=\sum_{k=0}^{N-1}\overline{z}^k\psi_k(z); \forall z\in\Omega_I$$ where $\varphi_k,\psi_k:\Omega_I\longrightarrow \C_I$ are holomorphic for all $k=0,...,N-1.$
By hypothesis, $\Omega$ is a slice domain, then it intersects the real line in a point $p_0\in\R$ such that there exists a ball of center $p_0$ in which we can expand the functions $\varphi_k$ and $\psi_k$ so that we have $$\displaystyle\varphi_k(z)=\sum_{j=0}^\infty (z-p_0)^j\alpha_{k_j}\text{ and }\psi_k(z)=\sum_{j=0}^\infty (z-p_0)^j\beta_{k_j}; \forall k=0,...,N-1,$$ where the coefficients $\alpha_{k_j}$ and $\beta_{k_j}$ do not depend on $I\in\Sq$.
Hence, $$\displaystyle f_I(z)=\sum_{k=0}^{N-1}\overline{z}^k\left(\sum_{j=0}^\infty(z-p_0)^j(\alpha_{k_j}+\beta_{k_j}J)\right).$$

Finally, the thesis follows from the arbitrariness of $I\in\mathbb{S}$.
\end{proof}
Therefore, we have the following characterization of slice polyanalytic functions on slice domains
\begin{cor}\label{carac3}
A function $f$ defined on a slice domain is slice polyanalytic of order $N$ if and only if it has the form \eqref{sum'}.
\end{cor}
\begin{proof}
This is a direct consequence of the Propositions \ref{carac1} and \ref{carac2}.
\end{proof}
The next results of slice regular functions that we shall extend to a higher order in this section are the counterparts of the identity principle, representation formula, extension lemma and the refined splitting lemma for slice polyanalytic functions.
\begin{thm}[Identity Principle]
Let $f$ and $g$ be two slice polyanalytic functions of order $N$ on a slice domain $\Omega\subset\Hq$. If, for some $I\in\Sq,$ $f$ and $g$ coincide on $U$ a subdomain of $\Omega_I$, then $f=g$ everywhere in $\Omega$.
\end{thm}
\begin{proof}
Note that $f$ and $g$ are slice polyanalytic functions of order $N$ on $\Omega$. Thus, we can write $$\displaystyle f(q)=\sum_{k=0}^{N-1}\overline{q}^kf_k(q) \text{ and } g(q)=\sum_{k=0}^{N-1}\overline{q}^kg_k(q); \forall q\in\Omega$$
where $(f_k)_{k=0,...,N-1}\text{ and } (g_k)_{k=0,...,N-1}$ are slice regular on $\Omega$. Note that thanks to the Splitting Lemma for slice polyanalytic functions we have that $f_I=F_1+F_2J$ and $g_I=G_1+G_2J$ where $J\in\mathbb{S}$ such that $I\perp J$ and $F_1,F_2,G_1,G_2$ are four $\C_I-$valued polyanalytic functions on $\Omega_I$. By hypothesis, we have  $f_I=g_I$ on $U$, so $F_1=G_1$ and $F_2=G_2$ on $U$ which is a subdomain of $\Omega_I$. Thus, from classical complex analysis we know that $F_1=G_1$ and $F_2=G_2$ everywhere on $\Omega_I.$ In particular, we get that $f_I=g_I$ everywhere on $\Omega_I$. Hence, $\overline{\partial}_I^{N-1}f_I=\overline{\partial}_I^{N-1}g_I$ on $\Omega_I$ which shows that $f_{N-1}$ coincides with $g_{N-1}$ on $\Omega_I$. However, $f_{N-1}$ and $g_{N-1}$ are slice regular. Then, making use of the Identity Principle for slice regular functions we have that $f_{N-1}=g_{N-1}$ everywhere on $\Omega$. Similarly, using the same arguments we show that $f_k=g_k$ on $\Omega$ for all $k=0,...,N-1.$ This ends the proof.
\end{proof}
 Inspired from the proof proposed in \cite{ColomboSabadiniStruppa2016} for slice regular functions, we can prove a representation formula for quaternionic slice polyanalytic functions:
\begin{thm}[Representation Formula] \label{RFSP}
Let $f$ be a slice polyanalytic function of order $N$ defined on an axially symmetric slice domain $\Omega\subset\Hq$. Let $J\in\Sq$, then for any $q=x+Iy\in\Omega$ the following equality holds :
$$
\displaystyle f(x+Iy)= \frac{1}{2}\left[f_J(x+Jy)+f_J(x-Jy)\right]+I\frac{J}{2}\left[f_J(x-Jy)-f_J(x+Jy)\right]
$$
Moreover, for all $x+yK\subset\Omega$, $K\in\Sq$, there exist two functions $\alpha,\beta$ independent of $I$, such that for any $K\in\Sq$ we have
$$\frac{1}{2}\left[f_K(x+yK)+f_K(x-yK)\right]=\alpha(x,y)$$
and $$\frac{1}{2}K\left[f_K(x-yK)-f_K(x+yK)\right]=\beta(x,y).$$
\end{thm}
\begin{proof}
Let $J$ be any imaginary unit in $\Sq$. For $q\in{\Omega}$, if $Im(q)=0$ the proof is obvious.  If not, we consider the function $\psi:\Omega\longrightarrow{\mathbb{H}}$, defined by \begin{center}
 $\psi(q):=\displaystyle\frac{1}{2}[f(Re(q)+J\vert{Im(q)}\vert)+f(Re(q)-J\vert{Im(q)}\vert)]
 +\displaystyle\frac{Im(q)}{\vert{Im(q)}\vert}\frac{J}{2}[f(Re(q)-J\vert{Im(q)}\vert)-f(Re(q)+J\vert{Im(q)}\vert)].$
\end{center}
We use the fact that $q=x+Iy$, $x,y\in\R, y\geq 0$ and $I=\displaystyle\frac{Im(q)}{\vert{Im(q)}\vert}$ and get  \begin{center}
$\psi(x+Iy):=\displaystyle\frac{1}{2}[f(x+Jy)+f(x-Jy)]+\displaystyle\frac{IJ}{2}[f(x-Jy)-f(x+Jy)]$
\end{center}
In particular, if $I=J$ we get \begin{center}
$\psi(x+Jy)=f_J(x+Jy)$
\end{center} this implies that $\psi\equiv{f_J}$ on $\Omega_J$ which is a domain by hypothesis. Hence, if we show that $\psi\in{\mathcal{SP}_N(\Omega)}$ we can use the identity principle for slice polyanalytic functions and conclude that $\psi$ and $f$ coincide everywhere on $\Omega$ and then
$$\displaystyle f(x+Iy)=\psi(x+Iy)=\frac{1}{2}[f(x+Jy)+f(x-Jy)]+\displaystyle\frac{IJ}{2}[f(x-Jy)-f(x+Jy)].$$ So, let us prove now that $\psi\in\mathcal{SP}_N(\Omega)$. Indeed, let $I\in\Sq$, note that the operators $\displaystyle\frac{\partial}{\partial x}$ and $\displaystyle I\frac{\partial}{\partial y}$ commute on $\mathcal{SP}(\Omega)$. Then, the binomial formula gives $$\displaystyle\overline{\partial}_I^N:=\frac{1}{2^N}\sum_{k=0}^N {N \choose k}I^k\frac{\partial^N}{\partial x^{N-k}\partial y^k}.$$ Therefore, by separating even and odd indices we obtain the following formula
$$\displaystyle\overline{\partial}_I^N=\frac{1}{2^N}\left[\sum_{k=0}^N (-1)^k{N \choose 2k}\frac{\partial^N}{\partial x^{N-2k}\partial y^{2k}}+I\sum_{k=0}^N (-1)^k{N \choose 2k+1}\frac{\partial^N}{\partial x^{N-2k-1}\partial y^{2k+1}}\right].$$

However, our aim is to show that $$\displaystyle\overline{\partial}^N_I\psi(x+Iy)=0.$$
To this end, let us set $$\displaystyle A_N(x,y)=\sum_{k=0}^N (-1)^k{N \choose 2k}\frac{\partial^N}{\partial x^{N-2k}\partial y^{2k}}\psi(x+Iy)$$ and $$\displaystyle B_N(x,y)=\sum_{k=0}^N (-1)^k{N \choose 2k+1}\frac{\partial^N}{\partial x^{N-2k-1}\partial y^{2k+1}}\psi(x+Iy).$$
Clearly, we have  \begin{equation}\label{dpsi}
\displaystyle\overline{\partial}^N_I\psi(x+Iy)=\frac{1}{2^N}(A_N(x,y)+IB_N(x,y))
\end{equation}
We keep in mind that the hypothesis $f$ slice polyanalytic of order $N$ gives the two following equations
\begin{equation}\label{a}
\begin{split}
  \displaystyle \sum_{k=0}^N (-1)^k{N \choose 2k}\frac{\partial^N}{\partial x^{N-2k}\partial y^{2k}}f(x+Jy)=\\ -J\sum_{k=0}^N (-1)^k{N \choose 2k+1}\frac{\partial^N}{\partial x^{N-2k-1}\partial y^{2k+1}}f(x+Jy)
\end{split}
\end{equation}
and
\begin{equation}\label{b}
\begin{split}
  \displaystyle \sum_{k=0}^N (-1)^k{N \choose 2k}\frac{\partial^N}{\partial x^{N-2k}\partial y^{2k}}f(x-Jy)=\\ J\sum_{k=0}^N (-1)^k{N \choose 2k+1}\frac{\partial^N}{\partial x^{N-2k-1}\partial y^{2k+1}}f(x-Jy).
\end{split}
\end{equation}

We have $$\displaystyle\psi(x+Iy)=\frac{1}{2}[f(x+Jy)+f(x-Jy)]+\displaystyle\frac{IJ}{2}[f(x-Jy)-f(x+Jy)].$$
Then, making use of the formulas \eqref{a} and \eqref{b} we obtain

\begin{equation}\label{AN}
\begin{aligned}
 \displaystyle A_N(x,y)={} & \frac{J}{2}\left[\sum_{k=0}^N (-1)^k{N \choose 2k+1}\frac{\partial^N}{\partial x^{N-2k-1}\partial y^{2k+1}} \left( f(x-Jy)-f(x+Jy)\right) \right]\\
      & -\frac{I}{2}\left[\sum_{k=0}^N (-1)^k{N \choose 2k+1}\frac{\partial^N}{\partial x^{N-2k-1}\partial y^{2k+1}} \left( f(x+Jy)+f(x-Jy)\right) \right]
\end{aligned}
\end{equation}
On the other hand, we have \begin{equation} \label{BN}
\begin{aligned}
 \displaystyle IB_N(x,y)={} & \frac{I}{2}\left[\sum_{k=0}^N (-1)^k{N \choose 2k+1}\frac{\partial^N}{\partial x^{N-2k-1}\partial y^{2k+1}} \left( f(x+Jy)+f(x-Jy)\right) \right]\\
&-\frac{J}{2}\left[\sum_{k=0}^N (-1)^k{N \choose 2k+1}\frac{\partial^N}{\partial x^{N-2k-1}\partial y^{2k+1}} \left( f(x-Jy)-f(x+Jy)\right) \right]
\end{aligned}
\end{equation}
Therefore, the formulas \eqref{AN} and \eqref{BN} combined with \eqref{dpsi} lead to $$\displaystyle\overline{\partial}^N_I\psi(x+Iy)=0.$$ This ends the proof.
\end{proof}
\begin{rem}
The proof of the second statement of Theorem \ref{RFSP} is similar to the one for slice regular functions which corresponds to $N=1$, see \cite{ColomboSabadiniStruppa2016}.
\end{rem}
Some immediate consequences of the representation formula for slice polyanalytic functions are the following :
\begin{cor}
Let $U\subset\Hq$ be an axially symmetric slice domain, $D\subset \R^2$ such that $x+yI\in U$ whenever $(x,y)\in D$ and let $f:U\longrightarrow \Hq$. Then, $f\in \mathcal{SP}_N(\Omega)$ if and only if there exist $\alpha,\beta:D\longrightarrow \Hq$ satisfying $\alpha(x,y)=\alpha(x,-y), \beta(x,y)=-\beta(x,-y)$ and $\overline{\partial_I}^N(\alpha+I\beta)=0$
such that $$f(x+yI)=\alpha(x,y)+I\beta(x,y).$$
\end{cor}
\begin{cor}
Let $U\subset\Hq$ be an axially symmetric slice domain and let $f:U\longrightarrow \Hq$ be a slice polyanalytic function. Then, for all $x,y\in \R$ such that $x+yI\in U$ there exist $a,b\in \Hq$ such that $$f(x+yI)=a+Ib$$ for all $I\in \mathbb{S}$.
\end{cor}
Another interesting fact that holds for slice polyanalytic functions is the Extension Lemma:
\begin{prop}[Extension] \label{Ext} Let $\Omega_I$ be a domain in $\C_I$ symmetric with respect to the real axis and such that $\Omega_I\cap\R\neq\emptyset$. If $f:\Omega_I\longrightarrow\Hq$ is polyanalytic of order $N$, then the function $$ext(f)(x+I_qy):=\displaystyle \frac{1}{2}\left[f(z)+f(\overline{z})\right]+I_q\frac{I}{2}\left[f(\overline{z})-f(z)\right]; z=x+Iy\in\Omega_I$$ is the unique slice polyanalytic extension of $f$ to the axially symmetric completion of $\Omega_I$ in $\Hq$. Moreover, if
$$f(z)=\displaystyle\sum_{k=0}^{N-1}\overline{z}^kh_k(z)$$ with $h_k:\Omega_I\longrightarrow\Hq$ are holomorphic functions. Then, we have
 $$ext(f)(q)=\displaystyle\sum_{k=0}^{N-1}\overline{q}^k ext(h_k)(q); \forall q=x+I_qy\in\Omega.$$
\end{prop}
\begin{proof}
Assume that $f$ is polyanalytic of order $N$ on $\Omega_I$. Then, we have $$f(z)=\displaystyle\sum_{k=0}^{N-1}\overline{z}^kh_k(z)$$ where $h_k:\Omega_I\longrightarrow\Hq$ are holomorphic functions. However, $\Omega_I$ is symmetric with respect to the real axis. Thus, according to the Extension Lemma for slice regular functions for any $k=0,...,N-1$ we can consider the slice regular functions defined by $$f_k(x+I_qy):=\displaystyle \frac{1}{2}\left[h_k(z)+h_k(\overline{z})\right]+I_q\frac{I}{2}\left[h_k(\overline{z})-h_k(z)\right]; z=x+Iy\in\Omega_I.$$
Let us consider $$g(x+I_qy)=\displaystyle\sum_{k=0}^{N-1}(x-I_qy)^kf_k(x+I_qy),$$ we shall prove that $$g(x+I_qy)=\displaystyle \frac{1}{2}\left[f(z)+f(\overline{z})\right]+I_q\frac{I}{2}\left[f(\overline{z})-f(z)\right]; z=x+Iy\in\Omega_I.$$
Indeed, first note that we have the two following equalities
\begin{equation}\label{i}
\displaystyle (x+I_qy)^k=\frac{1}{2}\left[(x+Iy)^k+(x-Iy)^k\right]+I_q\frac{I}{2}\left[(x-Iy)^k-(x+Iy)^k\right]
\end{equation}
and  \begin{equation}\label{ii}
\displaystyle (x-I_qy)^k=\frac{1}{2}\left[(x-Iy)^k+(x+Iy)^k\right]+I_q\frac{I}{2}\left[(x+Iy)^k-(x-Iy)^k\right].
\end{equation}
Then, by definition of $f_k$ we have $$g(x+I_qy)=\displaystyle \frac{C_N(x,y)+D_N(x,y)}{2}$$ where we have set $$C_N(x,y)=\displaystyle\sum_{k=0}^{N-1}(x-I_qy)^k\left[h_k(z)+h_k(\overline{z})\right]$$ and $$D_N(x,y)=\displaystyle\sum_{k=0}^{N-1}I_q(x-I_qy)^kI\left[h_k(\overline{z})-h_k(z)\right].$$
We replace $(x-I_qy)^k$ by its expression using the formula \eqref{ii} and get
\begin{equation}
\begin{split}
  \displaystyle C_N(x,y)=ext(f)(x+I_qy)+\frac{1}{2}\sum_{k=0}^{N-1}\left[z^kh_k(z)+\overline{z}^kh_k(\overline{z})\right]\\
+\frac{I_qI}{2}\sum_{k=0}^{N-1}\left[z^kh_k(z)-\overline{z}^kh_k(\overline{z})\right].
\end{split}
\end{equation}
On the other hand, after straightforward computations we obtain
\begin{equation}
\begin{split}
  \displaystyle D_N(x,y)=ext(f)(x+I_qy)-\frac{1}{2}\sum_{k=0}^{N-1}\left[z^kh_k(z)+\overline{z}^kh_k(\overline{z})\right]\\
+\frac{I_qI}{2}\sum_{k=0}^{N-1}\left[\overline{z}^kh_k(\overline{z})-z^kh_k(z)\right].
\end{split}
\end{equation}
Therefore, it follows that $$g(x+I_qy)=ext(f)(x+I_qy)$$ this ends the proof.
\end{proof}
Inspired from the book \cite{ColomboSabadiniStruppa2016}, we present the counterpart of the Refined Splitting Lemma for slice polyanalytic functions. First, let us consider the subclass of $\mathcal{SP}_N(\Omega)$ defined by $$\mathcal{N}_N(\Omega):=\lbrace{f\in{\mathcal{SP}_N(\Omega)} : f(\Omega\cap \C_I)\subset \C_I, \forall{I\in \mathbb{S}}}\rbrace.$$
Then, we have
\begin{prop}[Refined Splitting Lemma]
Let $\Omega$ be an axially symmetric slice domain in $\Hq$ and $f$ be a slice polyanalytic function of order $N$ on $\Omega$. Then, for any $I,J\in\mathbb{S}$ with $I\perp J$, there exist $\psi_\ell:\Omega_I\longrightarrow\C_I, \ell=0,..,3$ intrinsic polyanalytic such that: $$f_I(x+yI)=\psi_0(x+yI)+\psi_1(x+yI)I+\psi_2(x+yI)J+\psi_3(x+yI)K$$ where $K=IJ.$
\end{prop}
\begin{proof}
If $f$ is slice polyanalytic of order $N$, then we can write $$\displaystyle f(q)=\sum_{k=0}^{N-1}\overline{q}^kh_k(q) $$
with  $(h_k)_{k=0,..,N-1}$ are slice regular on $\Omega$. In particular, making use of the Refined Splitting Lemma for slice regular functions we have that for all $k=0,...,N-1:$ $$h_k(x+yI)=h_k^0(x+yI)+h_k^1(x+yI)I+h_k^2(x+yI)J+h_k^3(x+yI)IJ$$ where $h_k^\ell:\Omega_I\longrightarrow\C_I$ are holomorphic intrinsic functions for all $\ell=0,...,3.$
We have, $$f_I(z)=\displaystyle\sum_{k=0}^{N-1}\overline{z}^kh_k(z); \forall z\in\Omega_I.$$
Therefore, the thesis follows by considering the polyanalytic intrinsic functions defined by $$\psi_\ell(x+yI)=\displaystyle\sum_{k=0}^{N-1}(x-yI)^kh_k^\ell(x+yI); \forall \ell=0,...,3.$$
\end{proof}
As a consequence of the Refined Splitting Lemma, we have the following
\begin{thm}
Let $\Omega\subset \Hq$ be an axially symmetric slice domain and $\lbrace{1,I,J,IJ}\rbrace$ a basis of $\Hq$. Then, $$\mathcal{SP}_N(\Omega)=\mathcal{N}_N(\Omega)\oplus\mathcal{N}_N(\Omega)I \oplus\mathcal{N}_N(\Omega)J \oplus\mathcal{N}_N(\Omega) IJ.$$
\end{thm}
\begin{proof}
The Refined Splitting Lemma combined with the Extension Lemma for slice polyanalytic functions shows that we have
$$\mathcal{SP}_N(\Omega)=\mathcal{N}_N(\Omega)+\mathcal{N}_N(\Omega)I +\mathcal{N}_N(\Omega)J +\mathcal{N}_N(\Omega) IJ.$$
Moreover, we only need to use Proposition 2.7 in the book \cite{ColomboSabadiniStruppa2016} and the characterization of slice polyanalytic functions obtained in corollary \ref{carac3} to show that all the intersections between $\mathcal{N}_N(\Omega),\mathcal{N}_N(\Omega)I,\mathcal{N}_N(\Omega)J,\mathcal{N}_N(\Omega)IJ$ are reduced to zero. This ends the proof.
\end{proof}
As in the case of slice regular functions, namely $N=1$, we can introduce slice polyanalytic functions as a subclass of slice functions which are defined to be (see \cite{ColomboSabadiniStruppa2016}):
\begin{defn}
Let $U\subset\Hq$ be an axially symmetric open set. Functions of the form $f(q)=f(x+yI)=\alpha(x,y)+I\beta(x,y),$ where $\alpha,\beta$ are $\Hq-$valued functions such that $\alpha(x,-y)=\alpha(x,y), \beta(x,-y)=-\beta(x,y)$ for all $x+yI\in U$ are called slice functions.
\end{defn}
Then, we have the following:
\begin{defn}
An $N$ times differentiable slice function is said to be slice polyanalytic of order $N$ on an axially symmetric domain $\Omega$ if and only if for all $I\in\Sq,$ $f$ satisfies on $\Omega_I$ the equation $$\displaystyle \overline{\partial_I}^Nf(x+yI)=0.$$
\end{defn}
Inspired from the paper \cite{Altavilla2015}, we can show another version of the identity principle for slice polyanalytic functions without the hypothesis that the open set on which they are defined is a \textit{slice domain}. First, note that slice functions can be recovered by their values on two semi-slices, see the Representation Formula given by Theorem 2.4 in \cite{Altavilla2015}. We have the following
\begin{prop}Let $\Omega$ be an axially symmetric domain and let $f:\Omega\longrightarrow \Hq$ be a slice polyanalytic function. Assume that there exist $J,K\in\Sq,\text{ with } J\neq K$ and $U_J,U_K$ two subdomains of $\Omega_J^+$ and $\Omega_K^+$ respectively where $\Omega_J^+:=\Omega\cap\C_J^+$ and $\Omega_K^+:=\Omega\cap\C_K^+$. If $f=0$ on $U_J$ and $U_K$, then $f=0$ everywhere in $\Omega$.
\end{prop}
\begin{proof}
Let $f$ be a slice polyanalytic function on $\Omega$ such that $f=0$ on $U_J$ and $U_K$. Thus, since $U_J$ and $U_K$ are respectively subdomains of $\Omega_J^+$ and $\Omega^+_K$. It follows, from the Splitting Lemma for slice polyanalytic functions combined with the classical complex analysis that $f=0$ everywhere on $\Omega_J^+$ and $\Omega_K^+$. Then, we just need to use the Representation Formula which allows to recover a slice function  by its values on two semi-slices to complete the proof.
\end{proof}
\begin{rem}
This last remark on slice functions allows to define slice polyanalytic functions on axially symmetric domains which do not necessarily intersect the real line.
\end{rem}
\section{The quaternionic slice-polyanalytic Fock space}
In this section, we introduce the Fock space of slice polyanalytic functions on quaternions. Let $N\geq 1$ and $I\in\Sq$ we define the space
$$\displaystyle \mathcal{F}^{N}_{I}(\Hq):=\lbrace{f\in\mathcal{SP}_N(\Hq)/ \int_{\C_I}\vert{f_I(p)}\vert^2e^{-\vert{p}\vert^2}d\lambda_I(p)<\infty }\rbrace.$$
This space is endowed with the following inner product
$$\displaystyle\scal{f,g}_{\FNi}=\int_{\C_I}\overline{g_I(p)}f_I(p)e^{-\vert{p}\vert^2}d\lambda_I(p).$$
Then, we have the following:
\begin{prop} The set
$\mathcal{F}^{N}_{I}(\Hq)$ is a right quaternionic Hilbert space.
\end{prop}
\begin{proof}
The proof is based on the Splitting Lemma for slice polyanalytic functions, see Proposition \ref{split}. Indeed, let $(f_n)$ be a Cauchy sequence in $\FNi$. Choose $J\in\mathbb{S}$ such that $I\perp J$. Then, since $f_n$ are slice polyanalytic we have $f_{n,I}:=F_n+G_nJ \quad \forall n\in\N$ where $F_n$ and $G_n$ are polyanalytic functions on the slice $\C_I$ belonging to the classical polyanalytic Fock space $\mathcal{F}_N(\C_I)$. It is easy to see that $(F_n)_n$ and $(G_n)_n$ are Cauchy sequences in $\mathcal{F}_N(\C_I)$. Hence, there exists two functions $F$ and $G$ belonging to $\mathcal{F}_N(\C_I)$ such that the sequences $(F_n)_n$ and $(G_n)_n$ are converging respectively to $F$ and $G$. Let $f_I=F+GJ$ and consider $f=ext(f_I)$ we have then $f\in \FNi$ thanks to Proposition \ref{Ext}. Moreover, the sequence $(f_n)$ converges to $f$ with respect to the norm of $\FNi$. This ends the proof.
\end{proof}
\begin{prop} \label{independance}
Let $f\in\mathcal{SP}_N(\Hq)$ and $I,J\in\Sq$ two imaginary units. Then, we have the following $$ \frac{1}{2}\Vert{f}\Vert_{\mathcal{F}^{N}_{I}(\Hq)}\leq{\Vert{f}\Vert_{\mathcal{F}^{N}_{J}(\Hq)}}\leq 2\Vert{f}\Vert_{\mathcal{F}^{N}_{I}(\Hq)}.
  $$
\end{prop}
\begin{proof}
This is a consequence of the
Representation Formula, see Theorem \ref{RFSP}. Indeed, since $f$ is slice polyanalytic of order $N$ on $\Hq$ we have
  $$
\displaystyle f(x+Iy)= \frac{1}{2}\left[f(x+Jy)+f(x-Jy)\right]+I\frac{J}{2}\left[f(x-Jy)-f(x+Jy)\right].
$$
Then,
$$\left|f(x+Iy)\right| \leq  \left|f(x+Jy)\right|+ \left|f(x-Jy)\right|,$$
and therefore
\begin{align*}
\vert{f(x+Iy)}\vert^2 &\leq  \left(\vert{f(x+Jy)}\vert + \vert{f(x-Jy)}\vert\right)^2
\\& \leq 2\left(\vert{f(x+Jy)}\vert^2+\vert{f(x-Jy)}\vert^2\right)
\end{align*}
because $\left(\vert{f(x+Jy)}\vert - \vert{f(x-Jy)}\vert\right)^2 \geq 0$.
This implies that
\begin{align*}
\Vert{f}\Vert^2_{\FNi}
 \leq 2\left(\Vert{f}\Vert^2_{\FNj}+\Vert{f}\Vert^2_{\FNjj}\right).
\end{align*}
 However, since $\Vert{f}\Vert_{\FNj}=\Vert{f}\Vert_{\FNjj}$ we get $\Vert{f}\Vert^2_{\FNi}\leq 4\Vert{f}\Vert^2_{\FNj}$.
  By interchanging the roles of $I$ and $J$ we get also $\Vert{f}\Vert^2_{\FNi}\leq 4\Vert{f}\Vert^2_{\FNj}$. Finally, it follows that  $$\frac{1}{2}\Vert{f}\Vert_{\FNi}\leq{\Vert{f}\Vert_{\FNj}}\leq 2\Vert{f}\Vert_{\FNi}.$$
\end{proof}
\begin{cor}Given any $I,J\in\Sq,$ the slice polyanalytic Fock spaces $\FNi$ and $\FNj$ contain the same elements and have equivalent norms.
\end{cor}
\begin{rem} By the previous Corollary, the quaternionic slice polyanalytic Fock space is independent of the choice of the imaginary unit. Thus, we shall use the notation $\mathcal{F}^{N}_{Slice}(\Hq).$
\end{rem}
Let us fix $q\in\Hq$ and consider the evaluation mapping $$\Lambda_q:\FN\longrightarrow \Hq; f\mapsto\Lambda_q(f)=f(q).$$ Then, we have the following estimate on $\FN$ :
\begin{prop} \label{Growth}
Let $f\in\FN$ and $q\in\Hq$. Then, $$\vert{\Lambda_q(f)}\vert\leq \sqrt{N} e^{\frac{\vert{q}\vert^2}{2}}\Vert{f}\Vert_{\FN}.$$
\end{prop}
\begin{proof}
Let $I\in\Sq$ be such that $q\in\C_I$ and choose $J\in\Sq$ with $I\perp J$. Then, the Splitting Lemma yields $$f_I(z)=F(z)+G(z)J ;\text{ }\forall z\in\C_I$$ where $F$ and $G$  belong to $\mathcal{F}_N(\C_I)$. In particular, we have $$\vert{f(q)}\vert^2=\vert{F(q)}\vert^2+\vert{G(q)}\vert^2$$
However, we know from classical complex analysis that $$\vert{F(q)}\vert\leq \sqrt{N}e^{\frac{\vert{q}\vert^2}{2}}\Vert{F}\Vert_{\mathcal{F}_N(\C_I)}\text{ and } \vert{G(q)}\vert\leq \sqrt{N}e^{\frac{\vert{q}\vert^2}{2}}\Vert{G}\Vert_{\mathcal{F}_N(\C_I)}.$$
Therefore, $$\vert{f(q)}\vert\leq \sqrt{N}e^{\frac{\vert{q}\vert^2}{2}}\left(\Vert{F}\Vert^2_{\mathcal{F}_N(\C_I)}+\Vert{G}\Vert^2_{\mathcal{F}^N(\C_I)}\right)^{\frac{1}{2}}=\sqrt{N}e^{\frac{\vert{q}\vert^2}{2}}\Vert{f}\Vert_{\FN}.$$
\end{proof}

Proposition \ref{Growth} shows that all the evaluation mappings on $\FN$ are continuous. Then, the Riesz representation theorem for quaternionic right-linear Hilbert spaces, see \cite{BDS1982} shows that for any $q\in\Hq$ there exists a unique function $K_N^q\in\FN$ such that for any $f\in\FN$ we have $$f(q)=\scal{f,K^q_N}_{\FN}.$$
Let $J\in\Sq$ and $r\in\C_J$, then for $q=x+Iy$ and $z=x+Jy$ the corresponding reproducing kernel of the second kind is obtained by extending the kernel of the complex case. It is given by the following $$K_N:\Hq\times\Hq\longrightarrow \Hq$$
$$\displaystyle K_N(q,r):=\frac{1}{2}\left[K_N(z,r)+K_N(\overline{z},r)\right]+I\frac{J}{2}\left[K_N(\overline{z},r)-K_N(z,r)\right].$$

\begin{ex}
For $N=2$, we have
$$\displaystyle K_2(q,r)=ext[e^{z\overline{r}}(2-\vert{z-r}\vert^2)](q)$$
Then, we can check that $$K_2(q,r)=e_*(q\overline{r})(2-\vert{q-r}\vert^2).$$
\end{ex}This last example suggests to prove the following result
\begin{thm}\label{KerFock}
The set $\mathcal{F}^{N}_{Slice}(\Hq)$ is a right quaternionic reproducing kernel Hilbert space whose reproducing kernel is given by the following formula:

$$\displaystyle K_N(q,r)=e_*(q\overline{r})\sum_{k=0}^{N-1}(-1)^k{N \choose k+1}\frac{1}{k!}\vert{q-r}\vert^{2k}; \forall (q,r)\in\Hq\times\Hq.$$
\end{thm}
\begin{proof}
Fix $r\in\Hq$ such that $r$ belongs to the slice $\C_J$, we consider the function defined by $$\displaystyle F_N^r(q)=e_*(q\overline{r})\varphi_N(q,r) $$
 where  $$\varphi_N(q,r):=\sum_{k=0}^{N-1}(-1)^k{N \choose k+1}\frac{1}{k!}\vert{q-r}\vert^{2k}; \forall q\in\Hq.$$
Clearly $q\longmapsto e_*(q\overline{r})$ is slice regular on $\Hq$ with respect to the variable $q$. Moreover, we can check that  $\varphi_N(q,r)$ is a real valued slice polyanalytic function of order $N$ on $\Hq$ with respect to $q$. Thus, the observation in Proposition \ref{Mult} shows that $F_N^r$ is a slice polyanalytic function of order $N$ on $\Hq$ with respect to the variable $q$. Furthermore, note that the reproducing kernel of $\mathcal{F}^{N}_{Slice}(\Hq)$ is obtained by extending the classical one on the slice $\C_J$. This shows that $F_N^r(q)$ and $K_N(q,r)$ coincide on the slice $\C_J$ containing $r$. Hence, we have $K_N(q,r)=F_N^r(q)$ everywhere on $\Hq$ thanks to the Identity Principle for slice polyanalytic functions. This ends the proof.
\end{proof}
\begin{rem}
Note that the formula obtained in Theorem \ref{KerFock} could be written in terms of the generalized Laguerre polynomials $$\displaystyle L^\alpha_j(x):=\sum_{i=0}^j(-1)^i{j+\alpha \choose j-i}\frac{x^i}{i!}.$$
In particular, the reproducing kernel of the quaternionic slice polyanalytic Fock space $\FN$ will be written as $$K_N(q,r)=e_*(q\overline{r})L^1_{N-1}(\vert{q-r}\vert^2); \text{ } \forall (q,r)\in \Hq\times\Hq.$$
\end{rem}
\begin{rem}
For $N=1$, the space $\mathcal{F}^N_{Slice}(\Hq)$ is exactly the slice hyperholomorphic Fock space and the reproducing kernel obtained in Theorem \ref{KerFock} corresponds to the result obtained in \cite{AlpayColomboSabadini2014}.
\end{rem}
\section{The quaternionic slice polyanalytic Bergman space}
The slice polyanalytic Bergman space of the second kind on the quaternionic unit ball $\mathbb{B}$ is defined to be
$$\displaystyle \mathcal{A}^{N}_{Slice}(\mathbb{B}):=\lbrace{f\in\mathcal{SP}_N(\mathbb{B})/ \int_{\mathbb{B}_I}\vert{f_I(p)}\vert^2d\lambda_I(p)<\infty }\rbrace,$$
for $p=x+Iy$, $d\lambda_I(p)=dxdy$ is the usual Lebesgue measure on $\mathbb{B}_I=\mathbb{B}\cap\C_I$.
As we have seen in the previous section for the Fock space, we can use the same techniques involving the Splitting Lemma and Representation Formula for slice polyanalytic functions to prove that $\AN$ is a right quaternionic Hilbert space which does not depend on the choice of the slices. Furthermore,
for any $q\in\mathbb{B}$ and $f\in\AN$ we have the following estimate $$\vert{f(q)}\vert\leq \frac{N}{\sqrt{\pi}}\frac{\Vert{f}\Vert_{\AN}}{(1-\vert{q}\vert^2)}.$$
Hence, the Riesz representation theorem for quaternionic right-linear Hilbert spaces shows that $\AN$ has a reproducing kernel. The theory of quaternionic Bergman spaces of the second kind introduced in \cite{ColomboCervantesSabadini2015} suggests that the expression of the reproducing kernel of $\AN$ denoted by $B^N_S(q,r)$ is obtained making use of the extension operator. Indeed, let $r\in\mathbb{B}$ be fixed such that $r\in\mathbb{C}_J$, the expression of the kernel in the slice $\mathbb{B}_J$ is given in \cite{Balk1991} by
\begin{equation}
\displaystyle B_N^r(z)=\frac{N}{\pi(1-\overline{r}z)^{2N}}\sum_{k=0}^{N-1}(-1)^k{N \choose k+1}{N+k \choose N}\vert{1-\overline{r}z}\vert^{2(N-1-k)}\vert{z-r}\vert^{2k}.
\end{equation} Then, by definition, for any $q=x+Iy\in \mathbb{B}$ we have \begin{align*}
B_S^N(q,r) &= B_N^r(q)
\\&:= ext[B_N^r(z)](q).
\end{align*}
To give the explicit expression of $B_S^N(q,r)$, we consider first the function $f_N^r:\mathbb{B}_J\longrightarrow \C_J,$ depending on $r$ and defined by $$f_N^r(z)=\displaystyle \frac{N}{\pi}\frac{1}{(1-\overline{r}z)^{2N}}; \forall z\in\mathbb{B}_J.$$
We start by proving the following
\begin{lem}\label{Berlem}
For every fixed $r\in\mathbb{B}_J,$ the slice regular extension of $f_N^r(z)$ to the quaternionic unit ball $\mathbb{B}$ is given by
$$g_N^r(q)=P_N(q,r)Q_N(q,r);\text{ } \forall q\in\mathbb{B}$$ where
$$\displaystyle P_N(q,r)=\frac{N}{\pi}\sum_{k=0}^{2N}(-1)^k{2N \choose k}\overline{q}^k\overline{r}^k, \text{and } Q_N(q,r)=(1-2Re(q)\overline{r}+\vert{q}\vert^2\overline{r}^2)^{-2N}.$$
\end{lem}
\begin{proof}
Clearly, the function $f_N^r:z\mapsto f_N^r(z)$ is holomorphic on $\mathbb{B}_J$ for every fixed $r\in\mathbb{B}_J$. Then, by definition for $q=x+Iy \text{ and } z=x+Jy$ the slice regular extension of $f_N^r(z)$ to $\mathbb{B}$ is given by $$\displaystyle g_N^r(q)=\frac{1}{2}[f_N^r(z)+f_N^r(\bar{z})]+\frac{IJ}{2}[f_N^r(\bar{z})-f_N^r(z)].$$
We have \begin{align*}
\displaystyle\frac{f_N^r(z)+f_N^r(\bar{z})}{2} &=\frac{N}{2\pi}\left[\frac{1}{(1-\bar{r}z)^{2N}}+\frac{1}{(1-\bar{r}\bar{z})^{2N}}\right]
\\&= \frac{N}{\pi}\frac{\displaystyle\sum_{k=0}^{2N}(-1)^k{2N \choose k}\bar{r}^k \frac{(z^k+\bar{z}^k)}{2}}{(1-2Re(z)\bar{r}+\vert{z}\vert^2\bar{r}^2)^{2N}}
\\&= \frac{N}{\pi}\frac{\displaystyle\sum_{k=0}^{2N}(-1)^k{2N \choose k}\bar{r}^k Re(z^k)}{(1-2Re(z)\bar{r}+\vert{z}\vert^2\bar{r}^2)^{2N}}.
\end{align*}
Similarly, we obtain  $$
\displaystyle\frac{f_N^r(\bar{z})-f_N^r(z)}{2} =\frac{N}{\pi}\frac{J\displaystyle\sum_{k=0}^{2N}(-1)^k{2N \choose k}\bar{r}^k Im(z^k)}{(1-2Re(z)\bar{r}+\vert{z}\vert^2\bar{r}^2)^{2N}}.$$
Since $(1-2Re(z)\bar{r}+\vert{z}\vert^2\bar{r}^2)^{-2N}=Q_N(q,r)$, it follows by the formula of the extension operator that \begin{align*}
g_N^r(q)&= \displaystyle\frac{N}{\pi} \left[\sum
_{k=0}^{2N}(-1)^k{2N \choose k}(Re(z^k)-Im(z^k)I)\bar{r}^k\right]Q_N(q,r)
\\&=\displaystyle\frac{N}{\pi} \left[\sum
_{k=0}^{2N}(-1)^k{2N \choose k}\bar{q}^k\bar{r}^k\right]Q_N(q,r)
\\&=P_N(q,r)Q_N(q,r).
\end{align*}
This ends the proof.
\end{proof}
The next result gives the expression of the slice Bergman kernel of the second kind of the quaternionic slice polyanalytic Bergman space $\AN$:
\begin{thm}
The set $\AN$ is a right quaternionic reproducing kernel Hilbert space whose reproducing kernel is given by $$B_S^N(q,r)=P_N(q,r)Q_N(q,r)\psi_N(q,r); \text{ } \forall (q,r)\in \mathbb{B}\times\mathbb{B}$$
where
$$\displaystyle P_N(q,r)=\frac{N}{\pi}\sum_{k=0}^{2N}(-1)^k{2N \choose k}\overline{q}^k\overline{r}^k, \text{ } Q_N(q,r)=(1-2Re(q)\overline{r}+\vert{q}\vert^2\overline{r}^2)^{-2N}$$
and
$$\displaystyle\psi_N(q,r)=\sum_{k=0}^{N-1}(-1)^k{N \choose k+1}{N+k \choose N}\vert{1-\overline{r}q}\vert^{2(N-1-k)}\vert{q-r}\vert^{2k}.$$
\end{thm}
\begin{proof}
By definition of the second kind theory for quaternionic slice Bergman spaces, $B_S^N(q,r)$ is obtained making use of the extension operator. For any $r\in\mathbb{B}$ we consider the function $$h_N^r(q)=g_N^r(q)\psi_N(q,r) \text{ such that }g_N^r(q)=P_N(q,r)Q_N(q,r).$$ The polynomials $P_N(q,r)$ and $Q_N(q,r)$ are defined as in Lemma \ref{Berlem} and $$\displaystyle\psi_N(q,r)=\sum_{k=0}^{N-1}(-1)^k{N \choose k+1}{N+k \choose N}\vert{1-\overline{r}q}\vert^{2(N-1-k)}\vert{q-r}\vert^{2k}.$$ Then, clearly $g_N^r$ is slice regular on $\mathbb{B}$ with respect to the variable $q$ by construction according to Lemma \ref{Berlem}. We can see also that $\psi_N(q,r)$ is a real valued slice polyanalytic function of order $N$ on $\mathbb{B}$ with respect to $q$. Thus, the observation in Proposition \ref{Mult} shows that $h_N^r$ is a slice polyanalytic function of order $N$ on $\mathbb{B}$ with respect to $q$. Moreover, $h_N^r(q)$ and $B^N_S(q,r)$ coincide on the slice $\mathbb{B}_J$ containing $r$. Hence, thanks to the Identity Principle for slice polyanalytic functions $B_S^N(q,r)=h_N^r(q)$ everywhere on $\mathbb{B}$. This completes the proof.
\end{proof}
\begin{prop}
The kernel $B^N_S(q,r)$ can be written also in this second form $$B_S^N(q,r)=R_N(q,r)L_N(q,r)\psi_N(q,r); \text{ } \forall (q,r)\in \mathbb{B}\times \mathbb{B}$$
with
$$\displaystyle R_N(q,r)=(1-2qRe(r)+q^2\vert{r}\vert^2)^{-2N} \text{ and } L_N(q,r)=\frac{N}{\pi}\sum_{k=0}^{2N}(-1)^k{2N \choose k}q^kr^k.$$
\end{prop}
\begin{proof}
Set $\phi(q,r)=R_N(q,r)L_N(q,r)$ for all $q,r\in\mathbb{B}$. As a product of a rational function with real coefficients and a polynomial of order $2N$ with quaternionic coefficients on the right the function $\phi(.,r)$ is slice regular on $\mathbb{B}$ with respect to the variable $q$ for every $r\in\mathbb{B}$. Moreover, if $r\in\mathbb{B}$ is fixed on a slice $\C_J$ we can see that the restriction of $\phi(.,r)$ on $\mathbb{B}_J$ coincides with the function $f_N^r(z)=\displaystyle \frac{N}{\pi}\frac{1}{(1-\overline{r}z)^{2N}}$. Then, the Identity Principle for slice regular functions gives $$ext(f_N^r)(q)=R_N(q,r)L_N(q,r) \text{ for all } q,r\in\mathbb{B}.$$
The last equation leads to the desired result.
\end{proof}
\begin{rem}
For the particular case $N=1$, the results obtained in this section coincide with the results of \cite{ColomboCervantesSabadini2015} concerning the theory of the second kind for the slice hyperholomorphic Bergman spaces.
\end{rem}
\section{Concluding remarks}
We finish this paper by a brief discussion related to further developments of the theory of slice polyanalytic functions. First, we note that the pointwise multiplication of two slice polyanalytic functions is not slice polyanalytic, in general. In fact this fact appears also for $N=1$, namely for the case of slice regular functions. However, the $*$-product preserves the structure of slice regular functions. For slice polyanalaytic functions of the same order we can consider also a natural product denoted $*_N$ in order to preserve the structure. Indeed, let $f\text{ and } g$ be two slice polyanalytic functions of order $N$ on $\Omega$ such that $$\displaystyle f(q)=\sum_{k=0}^{N-1}\overline{q}^kf_k(q) \text{ and } g(q)=\sum_{k=0}^{N-1}\overline{q}^kg_k(q),$$ where $f_k,g_k$ are slice regular for all $k=0,..,N-1$. Then, we define $$\displaystyle f*_N g(q):=\sum_{k=0}^{N-1}\overline{q}^k(f_k*g_k)(q)$$ where $f_k*g_k$ stands for the classical $*$-product of slice regular functions. It turns out that the set $(\mathcal{SP}_N(\Omega),+,*_N)$ is a ring, so we wish to study further properties of this product in future researches.

 Furthermore, in the recent paper \cite{KahlerKuQian2017}, the authors introduced and studied the poly-Hardy space on the unit ball in the monogenic setting.  A natural problem would be to study the counterpart of the poly-Hardy space in this new slice polyanalytic setting. However, like in the classical complex case, this space would be trivial seen as subspace of $L^2(\mathbb{B})$.
\\ \\
\noindent{\bf Acknowledgements} \\
Daniel Alpay thanks the Foster G. and Mary McGaw Professorship in
Mathematical Sciences, which supported this research. Kamal Diki acknowleges the support of the project INdAM Doctoral Programme in Mathematics and/or Applications Cofunded by Marie Sklodowska-Curie Actions, acronym: INdAM-DP-COFUND-2015, grant number: 713485.

\newpage
\begin{center}
\textbf{CORRIGENDUM TO "ON SLICE POLYANALYTIC FUNCTIONS OF A QUATERNIONIC VARIABLE"}
\end{center}
In [2] we initiated the theory of polyanalytic functions of a quaternionic variable in the slice hyperholomorphic setting. While developing the associated functional calculus, see [1], and as it is discussed in [3, Remark 2.14], it turned out that  Definition 3.1 in [2] is too general to obtain all the expected results valid in the slice hypercomplex setting.   Thus,  the  slice polyanalytic functions of a quaternionic variable (or of a paravector variable, in the case of Clifford algebra-valued functions) have to be considered as a subclass of slice functions, see Definition 3.17 of [2].

The results presented in [2] for slice polyanalytic functions of a quaternionic variable stay valid when we use, instead of Definition 3.1 in [2],  Definition 3.17 of the same paper.
\begin{defn}[Slice  polyanalytic functions]\label{polycorrected}
Let $M\in \mathbb{N}$ and denote by $\mathcal{C}^M(U)$ the set of
continuously differentiable functions with all their derivatives up to order $M$ on an
axially symmetric open set $U\subseteq\mathbb{H}$.
We let $\mathcal{U} = \{ (x,y)\in\mathbb{R}^2\ :\ x+ I y\subset U\}$.
A function $F:U\to \mathbb{H}$
 is called a left
 slice function, if it is of the form
 \[
 F(q) = \alpha(x,y) + I\beta(x,y)\qquad \text{for } q = x + I y\in U
 \]
with the two functions $\alpha, \beta: \mathcal{U}\to \mathbb{H}$ that satisfy the compatibility conditions
$\alpha(x,-y) = \alpha(x,y)$, $\beta(x,-y) = -\beta(x,y)$.
If in addition $\alpha$ and $\beta$ are in $\mathcal{C}^M(U)$ and satisfy the poly Cauchy-Riemann equations of order $M\in \mathbb{N}$
 \begin{align}\label{CR}
\frac{1}{2^M}(\partial_x+I\partial_y)^M(\alpha(x,y) + I\beta(x,y))=0,\ \ \ {\rm for\ all} \ \ I\in \mathbb{S}
\end{align}
 then $F$ is called left slice polyanalytic function  of order $M\in \mathbb{N}$.
\end{defn}
The definition is easily adapted in the case of right slice polyanalytic functions. Note that a slice regular functions are obtained from the previous definition for $M=1$.
\\
Using the Definition \ref{polycorrected}, all the results of slice polyanalytic functions that were investigated in [2] are valid in the quaternionic setting, possibly with some changes in the proofs and in the list below we discuss all of them in detail. Note that the open sets on which the slice polyanalytic functions are defined are necessarily axially symmetric.
\\
(1) A slice function is said to be intrinsic when $\alpha$, $\beta$ in Definition \ref{polycorrected} are real-valued. Two slice functions can be multiplied with the $*$-product, see [4]. Thus Proposition 3.3 and its proof are valid.
\\
(2) Proposition 3.4, namely the Splitting Lemma, holds since its proof depends only on the splitting of the values of the functions. As a consequence also Proposition 3.5 holds since $f$ is clearly a slice function.
\\
(3)  The statement of Proposition 3.6 is valid but the proof follows the arguments used in the slice poly  monogenic setting, see Theorem 3.7 of [1]. We point out that Theorem 2.16 in [3] also establishes the result.
   \\
   The proof is as follows: let
 $f(q)=f(x+Iy)=\alpha(x,y)+I\beta(x,y)$ be a left slice polyanalytic of order $N$.
  By fixing a basis $1,e_1,e_2, e_1e_2$ of $\mathbb{H}$, and setting $e_0=1$, $e_3=e_1e_2$, we have
$\alpha=\sum_{\ell=0}^3 \alpha_{\ell} e_\ell$, $\beta=\sum_{\ell=0}^3 \beta_{\ell} e_\ell$, where the functions $\alpha_\ell, \beta_\ell$
are real-valued and are, respectively, even and odd in the second variable. Since the basis elements
$e_\ell$ are  linearly independent, the system expressing the slice polyanalyticity can be rewritten in terms of the
 components of $f$, in other words,
 each $\mathbb{C}_I$-valued function $F_\ell=\alpha_\ell+I\beta_\ell$ is polyanalytic and $\overline{F_{\ell}(x-Iy)}=F_{\ell}(x+Iy)$. By
the classical result applied to each function $F_\ell$, we have
$
F_{\ell}(x+Iy)=\sum_{k=0}^{N-1} (x-Iy)^k
f_{k,\ell}(x+Iy)
$
where the functions $f_{k,\ell}$ are $\mathbb{C}_{I}$-valued and
satisfy the Cauchy-Riemann system. Since $ \overline{F_{\ell}(x-Iy)}=\sum_{k=0}^{N-1} (x-Iy)^k
\overline{f_{k,\ell}(x-Iy)}=\sum_{k=0}^{N-1} (x-Iy)^k
f_{k,\ell}(x+Iy)$ we have $\overline{f_{k,\ell}(x-Iy)}=f_{k,\ell}(x+Iy)$ and so each $f_{k,\ell}$ is a slice function, as its real components satisfy the even-odd conditions.
We then deduce
\[
\begin{split}
f(x+Iy)&= \sum_{\ell=0}^3 F_\ell(x+Iy) e_\ell =
\sum_{\ell=0}^3\sum_{k=0}^{N-1} (x-Iy)^k f_{k,\ell}(x+Iy)e_\ell\\
&= \sum_{k=0}^{N-1} (x-Iy)^k f_{k}(x+Iy), \qquad f_k(x+Iy)= \sum_{\ell=0}^3 f_{k,\ell}(x+Iy)e_\ell,
\end{split}
\]
where the functions $f_k$ are
evidently left slice regular, and this concludes the proof.
\\
We note that Corollary 3.7 and its proof remains the same.
\\
(4) Theorem 3.8 in [2] remains valid, but in the proof the last five lines have to substituted by the observation that since $f$, $g$ are slice function, and since $f_I=g_I$, i.e. they coincide on $\Omega_I$, then they coincide on $\Omega$, by Proposition 5 in [4].
\\
(5) The representation formula in Theorem 3.9 is valid but there is no need of a proof since it is a consequence of the sliceness of a slice polyanalytic function, see [4]. Corollary 3.11 and 3.12 are correct but trivial.
\\
(6) The first part of Proposition 3.13 is not anymore needed so the statement become as below; the proof remains valid.
        \\
        {\em Let $\Omega_I$ be a domain in $\mathbb{C}_I$ symmetric with respect to the real axis and such that $\Omega_I\cap\mathbb{R}\neq\emptyset$. If
$$f(z)=\displaystyle\sum_{k=0}^{N-1}\overline{z}^kh_k(z)$$ with $h_k:\Omega_I\longrightarrow\mathbb{H}$ are holomorphic functions such that $\overline{h_k(\bar z)}=h_k(z)$. Then the unique slice polyanalytic extension of $f$ is
 $$ext(f)(q):=\displaystyle\sum_{k=0}^{N-1}\overline{q}^k ext(h_k)(q); \forall q=x+I_qy\in\Omega.$$}
\\
(7) In the sections 4, 5  of [2] concerning the function spaces, in order to have kernels that are slice functions we use  the $*$-product of (left) slice functions in the first variable, and we write the kernels in Theorems 4.6, Theorem 5.2 and Proposition 5.3 as
$K_N(q,r)=e_*(q\bar r)*\varphi_N(q,r)$ and $B_S^N(q,r)=P_N(q,r)Q_N(q,r)*\psi_N(q,r)=R_N(q,r)L_N(q,r)*\psi_N(q,r)$
where the slice functions (not necessarily real-valued) $\varphi_N, \psi_N$ are given by
$$
\varphi_N(q,r)=\sum_{k=0}^{N-1} (-1)^k {N\choose{k+1}}\frac{1}{k!}(\bar q q-q \bar r -\bar q r + \bar r r)^{k*},
$$
$$
\psi_N(q,r)=\sum_{k=0}^{N-1} (-1)^k {N\choose{k+1}} {{N+k}\choose{N}}(1-\bar q r-q\bar r+\bar q q  \bar r r)^{(N-1-k)*}*(\bar q q-q \bar r -\bar q r + \bar r r)^{k*}.
$$

{\bf References}

[1] D.~Alpay, F.~Colombo, K.~Diki, and I.~Sabadini.
\newblock {Poly slice monogenic functions, Cauchy formulas and the $PS$-functional calculus}.
\newblock{ ArXiv:2011.13912}, 2020.

[2] D.~Alpay, K.~Diki, and I.~Sabadini.
\newblock On slice polyanalytic functions of a quaternionic variable.
\newblock{Results in Mathematics}, 74(1):17, 2019.

[3] R. Ghiloni.
\newblock Slice-by-slice and global smoothness of slice regular and polyanalytic functions
\newblock { ArXiv:2011.09919}, 2020.

[4] R. Ghiloni, A. Perotti,
 {\em  Slice regular functions on real alternative algebras},
  Adv. Math., { 226} (2011),  1662--1691.





\end{document}